\documentclass[12pt]{amsart}

\usepackage{latexsym}
\usepackage[dvips]{graphics,epsfig}

\newtheorem{theorem}{Theorem}
\newtheorem{lemma}{Lemma}
\newtheorem{proposition}{Proposition}

\newtheorem{corollary}{Corollary}

\newcommand{\q}{\quad}
\newcommand{\qq}{\quad\quad}
\newcommand{\qqq}{\quad\quad\quad}
\newcommand{\qqqq}{\quad\quad\quad\quad}

\newcommand{\norm}[2]{{\left\| #1 \right\|}_{#2}}

\newcommand{\pnf}{+\infty}

\newcommand{\nf}{\infty}

\newcommand{\al}{\alpha}
\newcommand{\be}{\beta}
\newcommand{\ga}{\gamma}

\newcommand{\de}{\delta}
\newcommand{\ben}{\beta_n}

\newcommand{\ve}{\varepsilon}

\newcommand{\si}{\sigma}

\newcommand{\Om}{\Omega}

\newcommand{\rn}{{\mathbf R}^n}
\newcommand{\rone}{\mathbf R^1}
\newcommand{\rtwo}{\mathbf R^2}

\newcommand{\zz}{\mathbf Z}

\newcommand{\hone}{\mathbf H^1}

\newcommand{\sn}{\mathbf S^{n-1}}
\newcommand{\sone}{\mathbf S^1}

\newcommand{\lp}{L^{p}}

\newcommand{\lprn}{L^{p}(\rn)}

\newcommand{\lo}{L^{1}}
\newcommand{\lt}{L^{2}}

\newcommand{\cs}{\mathcal S}
\newcommand{\cm}{\mathcal M}

\newcommand{\pv}{\textup{p.v.}\,}

\newcommand{\intl}{\int\limits}

\newcommand{\liml}{\lim\limits}
\newcommand{\suml}{\sum\limits}
\newcommand{\supl}{\sup\limits}
\newcommand{\f}{\frac}
\newcommand{\df}{\displaystyle\frac}

\newcommand{\abssig}{\widehat{|\si_0|}}

\newcommand{\wt}{\widetilde}

\newcommand{\lbl}{\label}
\newcommand{\beq}{\begin{equation}}
\newcommand{\eeq}{\end{equation}}
\newcommand{\beqna}{\begin{eqnarray*}}
\newcommand{\eeqna}{\end{eqnarray*}}
\newcommand{\beqn}{\begin{equation*}}
\newcommand{\eeqn}{\end{equation*}}
\newcommand{\bp}{\begin{proof}}
\newcommand{\ep}{\end{proof}}
\newcommand{\bprop}{\begin{proposition}}
\newcommand{\eprop}{\end{proposition}}
\newcommand{\bt}{\begin{theorem}}
\newcommand{\et}{\end{theorem}}
\newcommand{\bex}{\begin{Example}}
\newcommand{\eex}{\end{Example}}
\newcommand{\bc}{\begin{corollary}}
\newcommand{\ec}{\end{corollary}}
\newcommand{\bl}{\begin{lemma}}
\newcommand{\el}{\end{lemma}}

\begin{document}

\title
[$L^p$ Bounds for  singular integrals and maximal singular integrals]
{$L^p$ bounds for  singular integrals \\
and maximal singular integrals with rough kernels}

\author{Loukas Grafakos}

\address{
Loukas Grafakos\\
Department of Mathematics\\
University of Missouri\\
Columbia, MO 65211, USA}

\curraddr{MSRI, Berkeley, CA 94720-5070}

\email{loukas@math.missouri.edu} 

\author{Atanas Stefanov}
\address{
Atanas Stefanov\\
Department of Mathematics\\
University of Missouri\\
Columbia, MO 65211, USA}

\curraddr{MSRI, Berkeley, CA 94720-5070}

\email{astefanov@pascal.math.missouri.edu}

\thanks{Grafakos's research at MSRI partially
supported by the NSF under grant DMS 9623120}




\begin{abstract}
Convolution type Calder\'on-Zygmund singular integral
operators with rough kernels $\pv \Om(x)/|x|^n$ are studied.
A condition on $\Om$ implying 
that the  corresponding  singular integrals and
maximal singular integrals map
$L^p \to L^p$ for $1<p<\nf$ is obtained. This condition is shown to be
different from  the condition $\Om\in H^1(\sn)$.
\end{abstract}

\maketitle

\section{Introduction and statements of results}

In this paper, $\Om$ will be a complex-valued
integrable function
defined on the sphere $\sn$, with mean value zero
with respect to surface measure.
Denote by $T_\Om$ the Calder\'on-Zygmund
singular integral operator  defined as follows:
\beq \lbl{1}
(T_\Om f)(x)=\liml_{\ve \to 0} \intl_{|y|>\ve}
\f{\Om (y/|y|)}{|y|^n}f(x-y)\, dy = \pv  \int_{\rn}
\f{\Om (y/|y|)}{|y|^n}f(x-y)\, dy,
\end{equation}
for $f$ in the Schwartz class $\cs (\rn)$. The limit in (\ref{1})
is  easily shown to exist for any $f$ a $C^1$ function on $\rn$
with some decay at infinity. 

For $\ve>0$,   denote by 
$$
(T^{\ve}_\Om f)(x)=
\int_{|y|>\ve}
\f{\Om (y/|y|)}{|y|^n}f(x-y)\, dy 
$$
the truncated singular integral associated with $T_\Om$ and by
$$
(T_\Om^* f)(x)=\sup_{\ve>0}|T^{\ve}_\Om f(x)|
$$
the maximal singular integral operator corresponding to this 
 $\Om$. 

Establishing the a priori  bound  $\|T_\Om^{\ve} f\|_{\lp}\le C \|f\|_{\lp}$
independently of  $f\in\cs(\rn)$ and of $\ve >0$,
leads to a (unique) extension of  $T^{\ve}_{\Om}$ on $ \lprn$.
Now, for $f\in \lprn$,  $T^{\ve}_{\Om}f$ converges in $\lp$  
as $\ve \to 0$ to some  $T_\Om f$ (which extends $T_\Om f$ defined 
in (\ref{1}) for $f\in \cs(\rn)$), 
and by Fatou's lemma, $T_\Om$ is a bounded operator on $\lp$.

A similar a priori bound for $T_\Om^*$ implies that for $f\in \lprn$, 
$T_{\ve}f$ converges (to $T_\Om f$) almost everywhere as $\ve \to 0$.

We now discuss $\lp$ boundedness properties of these operators. 
It is well known that if $\Om$ has some smoothness, then both
$T_\Om$ and $T_\Om^*$ extend to bounded operators on $\lprn$ for all
$1<p<\nf$. See \cite{Steinold} for details.
In this paper we shall be concerned with $\Om$ rough.
The method of rotations introduced by
Calder\'on and Zygmund \cite{CZ2}
implies that  that $T_\Om$ and $T_\Om^*$ map $\lprn\to \lprn$ for any
$\Om$ odd in $L^1 (\sn)$.
The  situation for general $\Om$'s is significantly
more involved. Calder\'on and Zygmund \cite{CZ2}
proved that  if
\beq\lbl{llogl}
\int_{\sn}|\Om (\theta)| \,  \textup{ln}(2+|\Om (\theta)|)
\, d\theta < \infty ,
\end{equation}
then $T_{\Om}$ and $T_\Om^*$ are bounded
operators on
$L^p$ for $1<p<\infty$.

Some years later,  condition (\ref{llogl}) above 
was independently improved by Connett \cite{connett} and
 Ricci and Weiss  \cite{RW}
 who showed that if 
\begin{equation}
\lbl{h1}
\Om \in H^1(\sn),
\end{equation}
then $T_{\Omega}$ maps $L^p(\rn)$ into itself for $1<p<\nf$. 
$H^1(\sn)$ here denotes  the $1$-Hardy space on the unit sphere 
 in the sense of Coifman and Weiss \cite{CW}; (this paper
contains a proof of this result in dimension $n=2$). See also 
\cite{graste} for a simple proof of this result on $\rn$. 

The $H^1$ condition (\ref{h1}) is also sufficient to imply that 
$T_\Om^*$ is bounded on $\lp$ for $1<p<\nf$. For a proof of this 
fact we refer the reader to \cite{graste} and also to Fan and Pan
\cite{fanpan} who recently obtained this result independently for a
more general class of operators. 

The main purpose of this paper is to present alternative conditions 
that imply $\lp$ boundedness for $T_\Om$ and $T_\Om^*$.  
If we examine the proof giving the formula of the
Fourier transform of $\pv \Om(x)/|x|^n$ we  observe that
the mild assumption 
\beq\lbl{mycondition}
\sup_{\xi\in \sn}
\int_{\sn} |\Om (\theta)| \ln \f{1}{| \theta \cdot \xi |}
\, d\theta <+ \nf,
\end{equation}
suffices to imply that
$(\pv \Om(x)/|x|^n)\widehat{\,\,\,\, }\,\,$ is a bounded  function, which 
is equivalent to saying  
 that $T_\Om$ maps $L^2(\rn)$ into itself. It is unknown to us
whether condition (\ref{mycondition}) implies $L^p$
boundedness for some $p\neq 2$.

Motivated by (\ref{mycondition}) we consider the family of
conditions
\beq\lbl{1+alpha}
\sup_{\xi\in \sn}
\int_{\sn} |\Om (\theta)| \left( \ln \f{1}{| \theta \cdot \xi|}
\right)^{1+\al}
\, d\theta <\pnf .
\end{equation}
for $\al>0$. We can show that if $\Om$ satisfies 
condition (\ref{1+alpha}) for some $\alpha >0$, then 
$T_\Om$ maps $\lprn$ into itself for some $p\neq 2$.
More precisely, we have the following theorem: 

\bt
Let  $\al>0$.
Let $\Om$ be a function in $L^1(\sn )$ with mean value zero
which satisfies condition (\ref{1+alpha})
for some $\al>0$.
Then $T_\Om$ extends to a bounded operator from $\lprn$ into itself
for $(2+\al)/(1+\al) <p < 2+\al$.
\et
As a corollary we obtain that if $\Om$ satisfies 
condition (\ref{1+alpha})
for all $\al>0$, then it maps $\lprn$ into itself for all $1<p<\nf$.
Regarding $T_\Om^*$ we can prove the following:

\bt
Let  $\al>1$.
Let $\Om$ be a function in $L^1(\sn )$ with mean value zero
which satisfies condition (\ref{1+alpha})
for some $\al>0$.
Then $T^{*}_\Om$ extends to a bounded operator from $\lprn$ into itself
for $1+3/(1+2\al) <p <2(2+\al)/3$.
\et

We conclude that if $\Om$ satisfies condition (\ref{1+alpha})
for all $\al>0$, then  $T^{*}_\Om$ maps $\lp$ to $\lp$ for all
$1<p<\nf$. We  don't know whether the ranges of indices
in Theorems 1 and 2 are sharp. More fundamentally, we do not know 
an example of an $\Om\in L^1(\sn)$ such that $T_\Om$ maps 
$\lp\to \lp$ for some given $p=p_0\ge 2$ but not for some other $p_1>p_0$.

In section 5, we show  that condition (\ref{1+alpha})
for all $\al>0$ is indeed disjoint from the $H^1$ condition (\ref{h1}).

\section{Boundedness of singular integrals}

The  theme of the 
proof  of Theorem 1 is based on  ideas developed by J. Duoandikoet- xea 
and J.-L. Rubio de Francia \cite{javirubio} to treat several other 
operators of this sort.
Define
$$
\si_k(x)= \f{\Om (x)}{|x|^n} \chi_{2^k\le |x|\le 2^{k+1}}, \qq k \in \zz.
$$
Observe that $\widehat{\si_k}(\xi)= \widehat{\si_0}(2^k\xi)$.
We calculate $\widehat{\si_0}(\xi)$. Set $\xi' =\xi/|\xi |$.
Expressing $\widehat{\sigma_0}$ in polar coordinates, we obtain 
\beq\lbl{5657}
\widehat{\si_0}(\xi)= \int_{\sn}\Om (\theta)
\left[ \int_1^2 e^{2\pi ir|\xi |
(\xi'\cdot \theta )}  \f{dr}{r} \right]\, d\theta.
\end{equation}
Using that $\Om$ has mean value zero, we deduce that
\beq\lbl{est1}
|\widehat{\si_0}(\xi)|\le  2\pi (\ln 2) \|\Om \|_{L^1} |\xi |=C|\xi |,
\end{equation}
which is a good estimate for $|\xi |\le 2$. For $|\xi |\ge 2$
observe the following: The integral inside brackets in
(\ref{5657}) is bounded by $\min \big(2, 3|\xi'\cdot \theta |^{-1}
|\xi |^{-1} \big)$. (Pick a $\theta$ so that $\xi'\cdot \theta \neq 0$.)
Therefore it must  satisfy the estimate
\begin{equation}\lbl{hnfiep}
\left| \int_1^2 e^{2\pi ir|\xi |
(\xi'\cdot \theta )} \f{dr}{r} \right| \le
\f{2 \left(\ln (\f32 |\xi'\cdot \theta |^{-1})  \right)^{1+\al} }
{( \ln |\xi | )^{1+\al} }.
\end{equation}
It follows from (\ref{hnfiep}) and (\ref{1+alpha}) that
\beq\lbl{est2}
|\widehat{\si_0}(\xi)|\le C (\ln |\xi |)^{-1-\al}\qq
\text{for $|\xi |\ge 2$.}
 \end{equation}
 Since $\si_k$ is obtained from
$\si_0$ by a suitable dilation, it follows  that  there exists a
constant $C>0$, such that  for all $k\in \zz$ the estimates below are valid:
\begin{align}\begin{split}\lbl{estk}
|\widehat{\si_k}(\xi)| &\le
C (\ln |2^k \xi |)^{-1-\al},\q
\text{for $2^k |\xi |\ge 2$, } \\
|\widehat{\si_k}(\xi)| &\le C 2^k |\xi |, \q
\qqqq\,\text{for $2^k |\xi |\le 2$.}
\end{split}\end{align}
Now let $\psi$ be a $C^\nf$ function supported in 
$\{x\in \rn :\,3/4\le |x|\le 9/4\}$ such that 
$\sum_{j\in \zz} (\psi (2^{j} \xi ))^2 =1$. Let $S_j$ be the
operator given on the Fourier transform by multiplication by
$\psi_j(\xi )= \psi (2^{j} \xi )$. Define
$$
T_jf= \sum_{k\in \zz} S_{j+k}(\si_k * S_{j+k}f).
$$
It is easy to see that the identity
$$
T_\Om f= \sum_{j \in \zz} T_jf
$$
is valid at least for $f$ in the Schwartz class.
Using a Fourier transform calculation, (\ref{estk}), and the
fact that $\psi_{j+k}$ is supported near the annulus
$|\xi | \sim 2^{-j-k}$, we obtain that 
$T_j$ are bounded on $\lt (\rn)$ with bound
$C2^{-j}$ for $j\ge 0 $ and $C(|j|)^{-1-\al}$ for $j\le -1 $.
In short
\beq\lbl{ltest}
\|T_jf\|_{\lt}
\le C(1+ |j|)^{-1-\al}\|f\|_{\lt} \qq\text{for all $j\in \zz$.}
\end{equation}
We will also need estimates for the following maximal operator
$$
f\to \si^*(f)= \sup_{k\in \zz} (|\si_k |*|f|).
$$
Without loss of generality we can assume 
that  
$\norm{\Om}{\lo(\sn)}=1$. It follows  that $\abssig(0)=1$.
Introduce a radial function in the Schwartz class $\Phi$, such that 
$\widehat{\Phi}(\xi)=1$ for $ |\xi|\le 2$ and $\widehat{\Phi}(\xi )=0$ 
for  $|\xi | >3$. Let us also introduce $\Phi_k$ defined by  
$\widehat{\Phi_k}(\xi)=\widehat{\Phi}(2^k \xi)$. Clearly we have 
\begin{equation}\label{maxfunc}
\si^*(f) \leq \sup_{k\in \zz}|(|\si_k |-\Phi_k)*|f||+
\sup_{k\in \zz}|\Phi_k*|f||. 
\end{equation}
Denote $\mu_k=|\si_k |-\Phi_k$. Since $\widehat{\mu_k}(0)=0$, the 
same proof giving (\ref{estk})  implies that 
\begin{align}\begin{split}
\lbl{estkp}
|\widehat{\mu_k}(\xi)| &\le C 2^k |\xi |, \q
\qqqq\,\text{for $2^k |\xi |\le 2$,}\\
|\widehat{\mu_k}(\xi)| &\le
C (\log |2^k \xi |)^{-1-\al},\q
\text{for $2^k |\xi |\ge 2$. } 
\end{split}\end{align}
Therefore we obtain from (\ref{maxfunc}) that 
\begin{equation}\label{mjd}
\si^*(f) \leq \sup_{k\in \zz}(\mu_k*|f|)+{\mathcal M}f \leq 
\left(\sum_k |\mu_k*|f||^2\right)^{1/2}+ {\mathcal M}f, 
\end{equation}
where ${\mathcal M}$ is the Hardy-Littlewood maximal function. 
Since for all $1<r<\nf$, 
\begin{equation}\lbl{13}
\Big\| \big(\sum_k (\mu_k*f)^2\big)^{1/2}\Big\|_{L^r}^r= \textup{Average}\q 
\big\| \sum_k \ve_k (\mu_k*f) \big\|_{L^r}^r, 
\end{equation}
over all choices of signs $\ve_k =\pm 1$, 
estimates for the square function on the right hand side 
of (\ref{mjd}) can be obtained from estimates on
integral operators of the form $g\to \sum_k \ve_k (\mu_k*g)$.
Now using  (\ref{estkp}) and (\ref{mjd})
we conclude  that $\si^*$ maps $\lt\to \lt$, whenever 
$\alpha>0$.
At this point we recall  the following lemma:

\bl\lbl{lemma} (See \cite{javirubio} p. 544)
If $\|\si^*(f)\|_{L^s}\le C
\|f\|_{L^s}$ and $\df{1}{2s}= \left| \df12- \df{1}{q}\right|$,
then for arbitrary  functions $g_k$ we have
\begin{equation*}
\big\|( \sum_{k \in \zz} |\si_k*g_k|^2)^{1/2} \big\|_{L^{q}} \le C
\big\|( \sum_{k \in \zz} |g_k|^2 )^{1/2} \big\|_{L^{q}}. 
\end{equation*}
\end{lemma}

Applying Lemma \ref{lemma} with $s=2$ and $ q=q_0=4$,
 we obtain that
\beq\lbl{lpest}
\|T_jf\|_{L^{q_0}} \le C
\big\| (\sum_{k\in \zz} |\si_k*S_{j+k}f|^2)^{1/2}
\big\|_{L^{q_0}}
\le C \big\| (\sum_{k\in \zz} |S_{j+k}f|^2 )^{1/2}
\big\|_{L^{q_0}} \le
C \|f\|_{L^{q_0}},
\end{equation}
where the middle inequality is a consequence of  Lemma
\ref{lemma} and the first and last inequalities follow from
the Littlewood-Paley theorem.

Interpolating between estimates (\ref{ltest}) and (\ref{lpest})
we obtain that
$$
\|T_jf\|_{L^p}\le C (1+|j|)^{-(1+\al)\theta_p}\|f\|_{\lp},
$$
where $1/p= \theta_p/2+(1-\theta_p)/q_0$.
Now observe that  $T_\Om=\sum_{j\in \zz} T_j$ maps $\lp\to \lp$
for all $p$'s for which  $p_1' < p< p_1$, where
$p_1= (4+4\al)/(2+\al)$ is the
unique solution of the equation $(1+\al)\theta_{p}=1$.
The same argument also gives that 
 $T_{\ve}f=\sum_k \ve_k (\mu_k*f)$ 
maps $\lp \to \lp$ for  $p_1' < p< p_1$   
uniformly on the choice of the signs $(\ve_j)$, $\ve_j=\pm 1.$ 
It follows that the square function in  (\ref{13}) is also bounded 
on $\lp$ for this range of $p$'s and hence so is 
$\si^*(f)$ by the estimate in (\ref{mjd}).
Thus we are in a position to apply  Lemma \ref{lemma} again
with $s$ in the interval $(p_1', p_1)$.

Now continue this way.
Fix $s_1\in (2, p_1)$ and let  $q_{1}$ be the
unique number bigger than $q_0=4$ which satisfies the equation
$1/2s_1'=|1/2-1/q_1|$. Apply Lemma \ref{lemma} with $s=s_1'$ and $q=q_1$.
As before we obtain that
 $T_\Om$
maps $\lp \to \lp$ for $p_2'< p< p_2$, where
$p_2$ is the unique solution of the equation
$(1+\al)\theta_{p}=1$, where  $\theta_p$ is given by
$1/p = \theta_p/2+(1-\theta_p)/q_1$
now. This bootstrapping argument leads to an
inductive definition of
three sequences  $2=p_0 < p_1< \dots \, $,
$2 < s_1<s_2< \dots \,$, and
 $4=q_0 < q_1< \dots \, $ such that for $k=1,2 , \dots $
$$
 p_{k-1} < s_{k}<p_{k}, \qq
\f{1}{p_{k}} - \f{1}{q_{k-1}}= \f{1}{1+\al} \left(
\f12-\f{1}{q_{k-1}}\right),\qq \f{1}{2s_k'}=\f12- \f{1}{q_{k}}.
$$
Let $b =\sup_k p_k$. The
above equations easily imply that $b = 2+\al$.
Therefore
$T_\Om$ maps $\lp$ to $\lp$ for $2\le p < 2+\al$.
The remaining range of $p$'s   follows by duality.

\section{Boundedness of maximal singular integrals}

We now prove Theorem 2. Below  we  use the 
same  notation  as in the previous section. 
Let 
\beqna
& & (T_k  f)(x)= 
 \int_{|y|>2^{k}}
\f{\Om (y)}{|y|^n}f(x-y)\, dy =\sum_{j=k}^{\infty}(\si_j *f)(x), \\
& & (T^*f)(x)=\sup_k|(T_kf)(x)|.
\eeqna
If $2^{k-1}\leq\ve<2^k$, then
$$
|(T_{\Om}^{\ve}f)(x)|\le|(T_k f)(x)|+\big|\int_{\ve<|y|<2^{k}}
\f{\Om (y)}{|y|^n}f(x-y)\, dy\big| \leq |(T_k f)(x)|+(|\si_k|*|f|)(x).
$$
From the proof of Theorem 1 we know that $\si^*$ maps $\lp\to\lp$ for 
$(2+\al)/(1+\al)<p<2+\al$. Since    
$$
|(T_{\Om}^*f)(x)|\leq |(T^*f)(x)|+\si^*(|f|)(x),
$$
it suffices to show that $T^*:\lp\to\lp$ for the claimed range of 
$p$'s, which is contained in the interval $((2+\al)/(1+\al),2+\al)$. 

With  $\Phi$ as in the previous section,  estimate 
\begin{equation}\lbl{RHS}
\sup_{k\in \zz}|(T_kf)(x)|\leq 
\sup_{k\in \zz}\left|\Phi_k*\sum_{j=k}^{\infty} \si_j*f\right|+
\sup_{k\in \zz}\left|(\delta-\Phi_k)*\sum_{j=k}^{\infty} \si_j*f\right|,
\end{equation}
where $\de $ is Dirac mass at the origin. 
It is  easy to see that
$$
\sup_{k\in \zz} \left|\Phi_k*\sum_{j=k}^{\infty} \si_j*f\right|\leq
 C \left(\cm(Tf)+\cm(f)\right), 
\qq \textup{(see \cite{javirubio}, p.548)}
$$
which implies $\lp$ bounds for  the first term on the right hand side
of (\ref{RHS}) for  \\ 
$(2+\al )/(1+\al)<p<2+\al$.
Control the second term   on the right hand side of
(\ref{RHS}) by 
$$
\sup_{k\in \zz} \left|(\delta-\Phi_k)*\sum_{j=0}^{\infty} \si_{j+k}*f\right|\leq
\sum_{j=0}^{\infty} Q_j(f),
$$
where 
$$
(Q_j f)(x)=  
\sup_{k\in \zz} \left|(\delta-\Phi_k)* \si_{j+k}*f\right| .
$$
To conclude the proof of Theorem 2, 
it suffices to show that for $j\ge 0$ we have 
\begin{align}
&\norm{Q_j f}{L^p}\le C\norm{f}{L^p},  \qqq 
       2\le p <2+\alpha,\lbl{130} \\
&\norm{Q_j f}{L^2}\le C (1+j)^{-\al} \norm{f}{L^2}.\lbl{131}
\end{align}
Then, a simple interpolation between (\ref{130}) 
and (\ref{131}) gives that $Q_j$ maps 
$\lp\to \lp$ with bound 
$C_\de (1+j)^{2\al (2+\al -\de-p)/p(\al-\de)}$, 
for any $\de>0$ small, and the conclusion of Theorem 2 
follows by summing on $j$. 

Now observe that
$$
|Q_j f| \le \sup_k |\si_{j+k}*f|+\sup_k|\Phi_k*\si_{j+k}*f|\le C(\si^*(f)+
\cm(\si^*(f))).
$$
Therefore $Q_j$ is bounded on $\lp$ whenever $\si^*$ is, that is 
$\norm{Q_jf}{\lp}\le C\norm{f}{\lp}$ when $2\le p<2+\al$ 
 and (\ref{130}) is proved.
To prove (\ref{131}) we need to exploit some orthogonality. We have 
$$
\|Q_j f\|_{\lt}^2\le \sum_k \|(\de-\Phi_k)*\si_{j+k}*f \|_{\lt}^2=
\textup{Average} 
\| \sum_k \ve_k  \big((\delta-\Phi_k)*\si_{j+k}*f \big)  \|^2_{\lt}, 
$$
where $\ve =(\ve_k)_{k}$ is a sequence of $\pm 1$'s. For a 
fixed sequence $\ve_k= \pm 1$, 
let us denote by
$$M_{j,k}f = \ve_k (\delta-\Phi_k)*\si_{j+k}*f.$$

We will need the following
\bl
Let $m\ge 1$, $j\ge 0$ and 
$ k_1\leq\ldots\leq k_{2m}$ be integers. Then 
$$
\norm{M_{j,k_1}\ldots M_{j,k_{2m}}}{2\to 2}\leq C^{2m} 
\prod_{i=1}^{2m}
\left(\f{1}{1+j+k_i-k_1}\right)^{1+\al} . 
$$
\el
\bp
Since $\widehat{\Phi_{k_1}}(\xi)$ vanishes for $2^{k_1}|\xi |\le 2$ 
we have,  
\beqna
\norm{M_{j,k_1}\ldots M_{j,k_{2m}}f}{L^2}^2&=&\int 
\prod_{i=1}^{2m}\left|1-\widehat{\Phi_{k_i}}(\xi)\right|^2
|\widehat{\si_{j+k_i}}(\xi)|^2|\widehat{f}(\xi)|^2 d\xi \\
&\leq& (C)^{2m}\int_{|\xi|\ge 2^{1-k_1}} \prod_{i=1}^{2m}\left[\f{1}
{\log\left(2^{j+k_i} |\xi |\right)}\right]^{2+2\al} |\widehat{f}(\xi)|^2 d\xi \\
&\leq& C^{2m}
 \prod_{i=1}^{2m}\left[\f{1}{1+j+k_i-k_1}\right]^{2+2\al} \norm{f}{L^2}^2,
\eeqna
where we used the first estimate in (\ref{estk}) in the 
last inequality above. 
\ep

Now we return to the proof of Theorem 2. We must  
show that
$\norm{M_j^{\ve,N}}{2\to 2}\leq C(1+j)^{-\al}$ uniformly on $N$ 
and $\ve=(\ve_k)$, where 
$$
M_j^{\ve,N}f=\sum_{k=-N}^{N}\ve_k M_{j,k}.
$$
Since $M_j^{\ve,N}$ are self adjoint operators, we have
\begin{gather*}
\norm{M_j^{\ve,N}}{2\to 2}^{2m}=\norm{(M_j^{\ve,N})^{2m}}{2\to 2}\leq
\sum_{-N\le k_1\le \ldots\leq k_{2m}\le N} 
\norm{M_{j,k_1}\ldots M_{j,k_{2m}}}{2\to 2}\\
\leq\sum_{-N\le k_1\le \ldots\leq k_{2m}\le N}C^{2m}\prod_{i=1}^{2m}
\left(\f{1}{1+j+k_i-k_1}\right)^{1+\al}\leq
\f{N C^{2m}}{(1+j)^{1+\al}}
\left(\f{1}{(1+j)^{\al}}\right)^{2m-1}\\
  \le N \f{C^{2m}}{1+j} (1+j)^{-2m\al}.
\end{gather*}
Taking $(2m)^{\textup{th}}$ roots and letting $m\to \nf$ we obtain 
$$
\norm{M_j^{\ve,N}}{2\to 2}\leq C (1+j)^{-\al}.
$$
This concludes the proof of (\ref{131}) and hence of 
Theorem 2.

\section{Examples}

It is easy to see that condition (\ref{1+alpha}) for all $\al>0$ 
contains the case $\Om \in L^q(\sn)$, $q>1$, considered by several authors, 
including \cite{javirubio}. However, it does not include the 
condition  $\Om \in L\textup{log}L(\sn)$ of Calder\'on and Zygmund. It is 
therefore natural to ask whether there exist examples of $\Om\notin
L\textup{log}L(\sn)$ which 
satisfy  (\ref{1+alpha}) for all $\al>0$. In this section we prove  
something more.

We construct an example to show that there exist 
integrable functions on $\sn$ with mean value zero which are 
not  in $H^1(\sn)$ but which 
satisfy  (\ref{1+alpha}) for all $\al>0$. Then we show that 
there exist functions that satisfy the converse. 

We begin with the converse which is easier. The function 
$$
\Om (\theta )= \sum_{k=2}^\infty \f{e^{ik\theta}}{(\log k)^2}
$$
belongs to $H^1(\sone)$ but it 
fails to satisfy condition (\ref{1+alpha}) for any $\al >0$.
Both assertions follow from the fact that $\Om (\theta )$ 
behaves like $\theta^{-1}
\log^{-2}(\theta^{-1})$ as $\theta\to 0+$ (See \cite{Zyg} p. 189).

We now construct an $\Om\in L^1(\sone) \setminus H^1(\sone)$
with mean value zero which  satisfies
condition (\ref{1+alpha}) for all $\al>0$.  
The example presented below is unavoidably complicated. The problem 
is that such a function must have an infinite number of spikes which 
are sufficiently far away from each other and which are (barely) 
integrable and have mean value zero. 

At this point we think of  $\sone$ as the interval 
$[0,1]$ via the identification 
\begin{equation}\lbl{ident}
\wt{\Om}(x) = \Om(\cos (2\pi x), \sin (2\pi x)),
\end{equation}
where $\Om$ is defined on $\sone$ and $\wt{\Om}$ on $[0,1]$.
It is not hard to see that  
under the identification given in (\ref{ident}), the condition 
$\Om\notin H^1(\sn)$  is equivalent to the 
fact that the Hilbert transform of $\wt{\Om}\chi_{[0,1]}$ 
is not in $L^1(\rone)$, and condition (\ref{1+alpha}) is equivalent to
\begin{equation}\lbl{1000}
\supl_{0\leq z\leq 1} \int_0^1|\wt{\Om}(x)| \ln^{1+\al}\f{1}{|x-z|} dx 
\leq C_\al<\infty.
\end{equation}
For a detailed justification of these facts see \cite{Stefanov}. 
Now let
\begin{align*}
 a_n &= (\ln n)^{-1} & b_n & =  e^{-\ga_n}\\
 \ga_n &=    e^{(\ln n)^{1/2}} & \de_n &=  e^{-\ga_n^{1/4}}\\
 d_n &= a_n+\de_n & c_n  &=  a_n-\de_n \\
 \be_n &= 1-(\ln n + \tfrac32 \ln \ga_n)\ga_n^{-1}&\q &\q 
 \end{align*}
Heuristically speaking,  
 $a_n$ is a sequence that decays slowly to zero, 
$c_n$ and $d_n$ are symmetric points about $a_n$ at 
distance $\de_n$, $(c_n-b_n,c_n)$ and $(d_n-b_n, d_n)$ are  
small intervals near $c_n$ and $d_n$ with length $b_n=e^{-\ga_n}$, 
where $(\ln n )^\ve << \ga_n<<n^\ve$ for all $\ve>0$, and the $\be_n$'s are 
powers that converge to one at a rate $\sim \ga_n^{-1}$. 
It is  easy to see that 
\begin{equation}\lbl{1008}
\f{b_n^{1-\be_n}}{1-\be_n}=\f{1}{n \ga_n^{1/2}(\ln  n +\f32 (\ln n )^{1/2})}
\sim \f{1}{n \ga_n^{1/2}\ln  n }, 
\end{equation} 
for $n$ large.  Now  let
$$
\wt{\Om}(x)=\suml_{n=10^9}^{\infty}\left(\f{1}{|x-c_n|^{\ben}}
\chi_{(c_n-b_n,c_n)}(x) - 
\f{1}{|x-d_n|^{\ben}}\chi_{(d_n-b_n,d_n)}(x)\right).
$$
We first verify  that condition (\ref{1000})  holds for all $\al>0$.
The worst possible $z$'s in (\ref{1000}) are the singularities of 
$\wt{\Om}$, i.e. the 
 points $z=c_n$, $d_n$, and $z=0$. By symmetry we consider only 
 $z=c_n$ and $z=0$. Fix $N\ge 10^9$ and consider 
 $z=c_N$.  We have 
$$
\int_0^1 |{\wt{\Om}}(x)| \ln^{1+\al}\f{1}{|x-c_N|} dx \le 
I_1(N)+I_2(N)+I_3(N)+I_4(N),
$$
where 
\begin{eqnarray*}
I_1 (N)& =&
 \suml_{n\neq N}
\intl_{c_n-b_n}^{c_n} \f{1}{|x-c_n|^{\be_n}} \ln^{1+\al}\f{1}{|x-c_N|} dx,  \\
I_2(N) &=&\suml_{n\neq N}\intl_{d_n-b_n}^{d_n} \f{1}{|x-d_n|^{\be_n}}
 \ln^{1+\al}\f{1}{|x-c_N|} dx, \\
I_3 (N)&=& \intl_{c_N-b_N}^{c_N} \f{1}{|x-c_N|^{\be_N}} 
\ln^{1+\al}\f{1}{|x-c_N|} dx ,\\
I_4 (N)&=& \intl_{d_N-b_N}^{d_N} \f{1}{|x-d_N|^{\be_N}} 
\ln^{1+\al}\f{1}{|x-c_N|} dx .
\end{eqnarray*}
Observe  that $I_2(N)\le CI_1(N)$ and that $I_4(N)\le I_3(N)$.  Also, it is easy to 
see that 
$$
\supl_{N\ge 10^9}I_3(N)
\leq C\supl_{N\ge 10^9} \f{b_N^{1-\be_N}}{1-\be_N}\ln^{1+\al}\f{1}{b_N}\leq 
C \supl_{N\ge 10^9}\f{\ga_N^{1+\al}}{N \ga_N^{1/2}\ln N}\leq C_\al.
$$
To control $\supl_{N\ge 10^9}I_1(N)$ 
we  need to show  that 
\begin{equation}\lbl{777777777}
\supl_{N\ge 10^9} \left[\suml_{n\neq N}\intl_{c_n-b_n}^{c_n}
 \f{1}{|x-c_n|^{\be_n}}
\ln^{1+\al}\f{1}{|x-c_N|} dx\right]\leq C_\al.
\end{equation}
Using that $|x-c_N|\sim|c_n-c_N|\sim |a_n-a_N|$ in the integrand above
 and (\ref{1008}), we conclude that (\ref{777777777}) will be a 
 consequence of 
\begin{equation}\lbl{1001}
\supl_{N\ge 10^9} \left[\suml_{n\neq N}\f{b_n^{1-\be_n}}{1-\be_n} 
\ln^{1+\al}\f{1}{|a_n-a_N|}\right]\leq C_\al.
\end{equation}
We have two cases. For $n > N, \q |a_n-a_N|\geq|a_{N+1}-a_N|\geq (N \ln^2
N)^{-1}$ and therefore
$$
\supl_{N\ge 10^9} \left[\suml_{n > N}\f{b_n^{1-\be_n}}{1-\be_n} 
\ln^{1+\al}\f{1}{|a_n-a_N|}\right]\leq	C 
\supl_{N\ge 10^9} \suml_{n > N}\f{\ln^{1+\al} (N \ln^2 N)}
{n \ga_n^{1/2}\ln n} \leq C_\al,
$$
the latter being an easy consequence of the integral test. For 
$10^9 \le n\le N-1$ we have 
\begin{align*}
&\suml_{n=10^9}^{ N-1}\f{b_n^{1-\be_n}}{1-\be_n} 
\ln^{1+\al}\!\f{1}{|a_n-a_N|}   \\ 
\le C &\suml_{n=10^9}^{ N-1}
\f{1}{n \ga_n^{1/2}\ln n}\ln^{1+\al}\!
\f{1}{|(\ln n)^{-1} - (\ln N)^{-1}|}=A(N)+B(N), 
\end{align*}
where $A(N)$ is the sum above of over the indices $10^9\le n< \ga_N $ and 
$B(N)$ is the sum over the the indices $\ga_N \le n\le N-1$. 
On $A(N)$ we have $|(\ln n )^{-1}- (\ln N)^{-1}|^{-1}
 \le C \ln n$, and thus $A(N)$ is 
clearly bounded independently of $N$. On $B(N)$ we have 
$|(\ln n)^{-1} - (\ln N)^{-1}|^{-1} \le C N (\ln N)^2$. Now estimate 
$\supl_{N\ge 10^9} B(N)$ by 
$$
C\supl_{N\ge 10^9}\ln^{1+\al}( N^2)\,\sum_{n\ge \ga_N} 
\f{1}{n\ga_n^{1/2}\ln n } \le 
C \supl_{N\ge 10^9}\f{\ln^{1+\al} (N^2)}{\ga_{\ga_N}^{1/3}} \le C,
$$
where we used the integral test to deduce the first inequality above. 
This concludes the proof of (\ref{1000}) when $z=c_N$. 
Condition (\ref{1000}) for $z=0$ is  is equivalent to 
the following inequality
$$
\suml_{n=10^9}^{\nf}\f{\ln^{1+\al}(\ln n)}{n \ga_n^{1/2}\ln n }\leq C_\al,
$$ 
which is certainly correct by the choice of our parameters. 
This proves that $\wt{\Om}$ satisfies condition (\ref{1000}) for all $\al>0$. 

We now prove  that $\wt{\Om}$ is not in the Hardy space 
$H^1$. Extend $\wt{\Om}$ 
to be equal to zero outside the interval $[0,1]$. 
Let $H$ be the usual Hilbert transform. 
Fix $N\ge 10^9$ and $y \in [d_N,d_N+b_N]$.
 Obviously 
\begin{equation}\lbl{6000}
\pi |(H\wt{\Om})(y)| \ge K_N(y)-L_N(y), 
\end{equation}
where 
\begin{eqnarray*}
K_N(y) &=& \left|\,\intl_{c_N-b_N}^{c_N}
\f{1}{|x-c_N|^{\be_N}}\f{1}{x-y}dx-
\intl_{d_N-b_N}^{d_N}\f{1}{|x-d_N|^{\be_N}}\f{1}{x-y}dx\right|, \lbl{10022} \\
L_N(y)&=& \suml_{n\neq N}\left|\, 
\intl_{c_n-b_n}^{c_n}\f{1}{|x-c_n|^{\be_n}}\f{1}{x-y}dx-
\intl_{d_n-b_n}^{d_n}\f{1}{|x-d_n|^{\be_n}}\f{1}{x-y}dx\right|.\lbl{1002}
\end{eqnarray*}
We first prove that 
\begin{equation}\lbl{6001}
\supl_{N\ge 10^9}\supl_{y \in [d_N,d_N+b_N]}L_N(y) \le C
\end{equation}
 Indeed,
\begin{eqnarray*}
& &\intl_{c_n-b_n}^{c_n}\f{1}{|x-c_n|^{\be_n}}\f{1}{x-y}dx=
\f{b_n^{1-\be_n}}{1-\be_n} \f{1}{(-b_n+c_n-y)} + \textup{smaller term} \\
& & \intl_{d_n-b_n}^{d_n}\f{1}{|x-d_n|^{\be_n}}\f{1}{x-y}dx=
\f{b_n^{1-\be_n}}{1-\be_n} \f{1}{(-b_n+d_n-y)} + \textup{smaller term}
\end{eqnarray*}
where the smaller terms are bounded by $ C b_n^{1-\be_n}/{(1-\be_n)}$ and 
$\sum\limits_{n\ge 10^9} b_n^{1-\be_n}/{(1-\be_n)}\leq C$.
Therefore
\begin{eqnarray*}
L_N (y)\le C \suml_{n\neq N}  
\f{b_n^{1-\be_n}}{1-\be_n}\f{|d_n-c_n|}{|a_n-a_N|^2}\leq 
C \suml_{n\neq N}\f{b_n^{1-\be_n}}{1-\be_n}\f{\de_n}{|a_n-a_N|^2}
\end{eqnarray*}
and thus it remains to prove that 
\begin{equation}\lbl{4000}
\supl_{N\ge 10^9}
\suml_{n\neq N}\f{\de_n}{n \ga_n^{1/2}\ln n } 
\f{1}{((\ln n)^{-1}-(\ln N)^{-1})^2} \leq C.
\end{equation}
The sum in (\ref{4000}) for $n>N$ is bounded by 
$$
\suml_{n> N}\f{\de_n}{n \ga_n^{1/2}\ln n} 
\f{1}{((\ln n)^{-1}-(\ln N)^{-1})^2}\leq
 N^2 \ln^4 N \suml_{n> N}\f{\de_n}{n \ga_n^{1/2}\ln n } \leq  C, 
$$
uniformly in $N\ge 10^9$. 
Split  the sum in (\ref{4000}) for  $n<N$ into the sum $A'(N)$ over  
the indices $10^9\le n < \ga_N$ and the sum $B'(N)$ over the 
indices $\ga_N\le n \le N-1$.
 Using that when $10^9\le n < \ga_N$    
 we have $|(\ln n)^{-1} - (\ln N)^{-1}|^{-1} \le C \ln n$ 
we conclude that $A'(N)$ is bounded independently of $N$. 
When $\ga_N\le n \le N-1$ we have 
$|(\ln n)^{-1} - (\ln N)^{-1}|^{-1} \le C N (\ln N)^2$ and hence 
$$
\supl_{N\ge 10^9} 
B'(N)\le C\supl_{N\ge 10^9} N^5 \,\sum_{n\ge \ln N} 
\f{1}{n\ga_n^{1/2}(\ln n ) e^{\ga_n^{1/4}}  }\le C,  
$$
which follows from the integral test. 
This proves (\ref{4000}) and hence  $L_N(y)$  is bounded uniformly in $N$. 

Now we turn our attention to $K_N(y)$. Observe
 that the following inequality holds 
$$
\intl_{d_N-b_N}^{d_N}\f{1}{|x-d_N|^{\be_N}}\f{1}{y-x}dx\geq
\f{3}{2}\intl_{c_N-b_N}^{c_N}\f{1}{|x-c_N|^{\be_N}}\f{1}{y-x}dx,
$$
because of the proximity of $y$ to the support of the first integral. 
Therefore
$$
|K_N(y)|\geq c \intl_{d_N-b_N}^{d_N}
\f{1}{|x-d_N|^{\be_N}}\f{1}{y-x}dx - C
$$
when $y \in [d_N,d_N+\de_N]$. Integrate over this set to  obtain 
\begin{gather}\begin{split}\lbl{2000}
\intl_{d_N}^{d_N+\de_N}| K_N(y)|dy \ge\q\q\q\q\q\q\q\q\q\q\q\q \\
c\left|\, \intl_{d_N-b_N}^{d_N}\f{1}{|x-d_N|^{\be_N}}\ln(d_N+\de_N-x)dx
-\intl_{d_N-b_N}^{d_N}\f{1}{|x-d_N|^{\be_N}}\ln(d_N-x)dx -\right| C\de_N.
\end{split}\end{gather}
We clearly have that 
\begin{equation}\lbl{2001}
\left|\, \intl_{d_N-b_N}^{d_N}\f{1}{|x-d_N|^{\be_N}}\ln(d_N+\de_N-x)dx\right|
\leq 
C |\ln \de_N | \f{b_N^{1-\be_N}}{1-\be_N}\leq 
\f{\ga_N^{1/4}}{N \ga_N^{1/2}\ln N},
\end{equation}
while the the crucial fact is that
\begin{equation}\lbl{2002}
\left|\, \intl_{d_N-b_N}^{d_N}\f{1}{|x-d_N|^{\be_N}}\ln(d_N-x)dx\right|
\geq C|\ln b_N|\f{b_N^{1-\be_N}}{1-\be_N}\geq\f{\ga_N^{1/2}}{N \ln N}.
\end{equation}
Combining (\ref{6000}), (\ref{6001}), (\ref{2000}), 
(\ref{2001}), and (\ref{2002}) 
we obtain 
\begin{align*} 
\norm{H\wt{\Om}}{L^1}\geq &\suml_{N\ge 10^9} \intl_{d_N}^{d_N+\de_N}
|(H\wt{\Om})(y)|dy  \\
\geq c &\suml_{N\ge 10^9}
\f{\ga_N^{1/2}}{N\ln N}- C\suml_{N\ge 10^9}
\f{1}{N \ga_N^{1/4}\ln N}-C\suml_{N\ge 10^9} \de_N=\infty.
\end{align*}
This proves that $\wt{\Om}\notin H^1([0,1])$.

\end{document}